\documentclass[11pt]{article}
\usepackage{amssymb}
\usepackage{graphicx}
\usepackage{amsmath}
\usepackage{xcolor}%color
\usepackage{colortbl,booktabs}%rule
\usepackage{multirow}
\usepackage{subfigure}
\usepackage{multicol}

\def\R{{\mathbb R}} 
\def\Z{{\mathbb Z}} \def\N{{\mathbb N}}
\def\Q{{\mathbb Q}}

\def\bb{\begin}
\def\bc{\begin{center}}
\def\ec{\end{center}}
\def\ba{\begin{array}}
\def\ea{\end{array}}

\def\be{\begin{equation}}     \def\ee{\end{equation}}
\def\bea{\begin{eqnarray}}    \def\eea{\end{eqnarray}}
\def\beaa{\begin{eqnarray*}}  \def\eeaa{\end{eqnarray*}}

\def\hh{\!\!\!\!}             \def\EM{\hh &   &\hh}
\def\EQ{\hh & = & \hh}        
      
\def\LT{\hh & < & \hh}        
\def\NE{\hh & \ne & \hh}      \def\AND#1{\hh & #1 & \hh}

              \def\vp{\varphi}

\def\vr{\varrho}

\def\nn{\nonumber}            
\def\oo{\infty}               \def\d{\cdot}
\def\dd{\cdots}               \def\pa{\partial}
\def\q{\quad}                 \def\qq{\qquad}
\def\f{\frac}                 
\def\z{\left}                 \def\y{\right}

\def\rd{\,{\rm d}}
\def\dt{\,{\rm d}t}
\def\ds{\,{\rm d}s}
\def\du{\,{\rm d}u}

\def\ifl{\iffalse}
\def\qqf{\qquad \forall}
\def\andq{\quad \mbox{ and } \quad}
\def\Proof{\noindent{\bf Proof} \quad}
\def\qed{\hfill $\Box$ \smallskip}

\def\lb{\label}
\def\x#1{(\ref{#1})}

\def\matt#1#2#3#4{\z( \ba{cc} #1 & #2 \\ #3 & #4 \ea \y)}
\def\pas#1#2{\left.#1\right|_{#2}}

\def\bu{$\bullet$\ }
\def\buu{$\bullet\bullet$\ }
\def\pp{{m,p}}
\def\vpm{\eta_\pp}
\def\T{{2m\pi}}
\def\imt{\int_0^\T}
\def\oh{\frac{\footnotesize 1}{\footnotesize 2}}
\def\oq{\frac{\footnotesize 1}{\footnotesize 4}}

\def\ul{\underline}
\def\cc{C}
\def\ss{S}
\def\X{\ul{X}}

\begin{document}

\begin{center}
{\LARGE On the Stability of Symmetric Periodic Orbits of \\
the Elliptic Sitnikov Problem}\\
\vskip 0.3cm

Xiuli Cen\footnote{This author is supported by the National Natural Science Foundation of China (Grant No. 11801582).}\\
School of Mathematics (Zhuhai), Sun Yat-sen University, \\
Zhuhai, Guangdong 519082, China \\
E-mail: {\tt cenxiuli2010@163.com} \\
\vskip 0.2cm

Xuhua Cheng\footnote{This author is supported by the National Natural Science Foundation of China (Grant No. 11601257).} \\
Department of Applied Mathematics, Hebei University of Technology, \\ Tianjin 300130, China\\
E-mail: {\tt chengxuhua88@163.com}
\vskip 0.2cm

Zaitang Huang\footnote{This author is supported by the Guangxi Natural Science Foundation (Grant No. 2018JJA110052).}\\
School of Mathematics and Statistics, Nanning Normal University,\\ Nanning 530023, China\\
E-mail: {\tt zaitanghuang@163.com}\\
\vskip 0.2cm

Meirong Zhang\footnote{Correspondence author. This author is supported by the National Natural Science Foundation of China (Grant No. 11790273).}\\
Department of Mathematical Sciences, Tsinghua University, \\ Beijing 100084, China\\
E-mail: {\tt zhangmr@tsinghua.edu.cn}
%%\today
\end{center}
\vskip 0.2cm

    \bb{abstract}
Motivated by the recent works on the stability of symmetric periodic orbits of the elliptic Sitnikov problem, for time-periodic Newtonian equations with symmetries, we will study symmetric periodic solutions which are emanated from nonconstant periodic solutions of autonomous equations. By using the theory of Hill's equations, we will first deduce in this paper a criterion for the linearized stability and instability of periodic solutions which are odd in time. Such a criterion is complementary to that for periodic solutions which are even in time, obtained recently by the present authors. Applying these criteria to the elliptic Sitnikov problem, we will prove in an analytical way that the odd $(2p,p)$-periodic solutions of the elliptic Sitnikov problem are hyperbolic and therefore are Lyapunov unstable when the eccentricity is small, while the corresponding even $(2p,p)$-periodic solutions are elliptic and linearized stable. These are the first analytical results on the stability of nonconstant periodic orbits of the elliptic Sitnikov problem.
    \end{abstract}



%%%\tableofcontents
%%%\newpage

\section{Introduction} \setcounter{section}{1} \setcounter{equation}{0} \lb{main-result}

The elliptic Sitnikov problem, denoted by $(S_e)$, is the simplest model in the restricted $3$-body problems \cite{S60}. By assuming that the two primaries with equal masses are moving in a circular or an elliptic orbit of the $2$-body problem of the eccentricity $e\in [0, 1)$, the Sitnikov problem describes the motion of the infinitesimal mass moving on the straight line orthogonal to the plane of motion of the primaries, %% which passes through their center of mass.
whose governing equation was given in \cite{BLO94, LO08} and will be stated as Eq. \x{se} in \S \ref{Sitnikov} of this paper. When $e=0$, $(S_0)$ is called the circular Sitnikov problem, whose equation, stated as Eq. \x{s0}, is an autonomous scalar Newtonian or Lagrangian equation. For $e\in(0,1)$, the equation for $(S_e)$ is a nonlinear scalar Newtonian equation which is $2\pi$-periodic in time.

There is a long history and a rigorous study on motions of problem $(S_e)$, covering the following topics.

\bu {\bf Oscillation and expressions of motions:} The motions of the circular Sitnikov problem can be expressed using various elliptic functions in an implicit way \cite{BLO94, F03, LS90, S60}. It is also found that the elliptic Sitnikov problem admits oscillatory motions. See the bibliography of \cite{LO08} for some historic references on this topic.

\bu {\bf Existence and construction of periodic orbits:} Due to the symmetries of the elliptic Sitnikov problem, many interesting periodic orbits have been obtained in \cite{BLO94, LO08, LS80, O16, OR10}, mainly by using the bifurcation method and global continuation.

\bu {\bf Stability and linearized stability of motions:} This is a central topic in dynamical systems \cite{O17, SM71}. For example, $(S_e)$ has the origin as an equilibrium which can be considered as a $2\pi$-periodic solution. In case the equilibrium is elliptic, its Lyapunov stability can be studied using the third order approximation developed by Ortega \cite{O96} and extended in \cite{LLYZ03}. See \cite[\S 6]{LO08} for details. As for nonconstant, even (in time) periodic solutions of $(S_e)$ which are emanated from the corresponding solutions of $(S_0)$, the stability and linearized stability are studied in very recent papers \cite{GNR18, GNRR18, M18, ZCC18}. Most of these are based on the theory for Hill's equations. Though some analytical formulas have been derived, many results of these are numerical due to the difficulties caused by nonconstant periodic solutions.

In this paper we continue the study for the stability and linearized stability of nonconstant, symmetric (in time) periodic solutions of $(S_e)$. Our aim is to provide some analytical results. In order to make such an analytical approach be applicable to more general problems, we consider the following second-order nonlinear scalar Newtonian equation
    \be \lb{xe}
    \ddot x+ F(x,t,e)=0.
    \ee
Here $F(x,t,e)$ is a smooth function of $(x,t,e)\in \R^3$ fulfilling the following symmetries
%%and the periodicity as in \x{Sy1}. Moreover, we assume that when $e=0$,
    \be \lb{Sy1}
    \z\{\ba{l} F(-x,t,e) \equiv -F(x,t,e), \\
    F(x,-t,e) \equiv F(x,t,e),\\
    F(x,t+2\pi,e) \equiv F(x,t,e),\\
    F(x,t,0)\equiv f(x), \\
    x f(x)> 0\mbox{ for }x\ne 0.
    \ea\y.
    \ee
These symmetries are verified by the Sitnikov problem $(S_e)$. In particular, when $e=0$, the starting equation
    \be \lb{x}
    \ddot x+f(x)=0
    \ee
is autonomous and has the unique equilibrium $x=0$. Obviously, $f(x)$ is also odd in $x$.

Let $m, \ p\in \N$ be integers. We say that $x(t)$ is an $(\pp)$-periodic solution of Eq. \x{xe}, if $x(t)$ is a $\T$-periodic solution of \x{xe} and has precisely $2p$ zeros in intervals $[t_0,t_0+\T)$, $t_0\in \R$.
%%Obviously, if $x(t)$ is $(\pp)$-periodic, then it is also $(m n, p n)$-periodic for any $n\in \N$.

Because of the autonomy and the complete integrability, all $(\pp)$-periodic solutions of Eq. \x{x} are clear. In particular, with suitable choice of $(\pp)$, Eq. \x{x} admits the $(\pp)$-periodic solutions $\vp_\pp(t)$ and $\phi_\pp(t)$, which are respectively even and odd in time $t$. These are the symmetric $(\pp)$-periodic solutions of Eq. \x{x} we are interested in. Due to the autonomy of Eq. \x{x}, both $\vp_\pp(t)$ and $\phi_\pp(t)$ have the minimal period $\T/p$.

From bifurcation theory, it is known that, under some non-degeneracy conditions, Eq. \x{xe} admits families of $(\pp)$-periodic solutions $\vp_\pp(t,e)$ and $\phi_\pp(t,e)$, $0 \le e \ll 1$, such that
    \[
    \z\{ \ba{l}
    \vp_\pp(t,0)\equiv\vp_\pp(t), \mbox{ and } \phi_\pp(t,0)\equiv\phi_\pp(t),\\
    \vp_\pp(t,e) \mbox{ is even in } t, \ \vp_\pp(0,e)>0, \mbox{ and } \vp_\pp(t+m\pi,e) \equiv - \vp_\pp(t,e), \\
    \phi_\pp(t,e) \mbox{ is odd in } t, \ \dot\phi_\pp(0,e)>0, \mbox{ and } \phi_\pp(t+m\pi,e) \equiv - \phi_\pp(t,e).
    \ea
    \y.
    \]
They are called the even and the odd $(\pp)$-periodic solutions of Eq. \x{xe}, respectively.
%%In fact, these solutions are $m\pi$-anti-periodic in $t$.
Generally speaking, when $e>0$, $\vp_\pp(t,e)$ and $\phi_\pp(t,e)$ have the minimal period $\T$, not $\T/p$. For more details, see Theorem \ref{M1}. For the elliptic Sitnikov problem $(S_e)$, such symmetric periodic solutions have been studied extensively in \cite{BLO94, LO08, O16}. Moreover, some interesting global continuations of these solutions are also obtained. See, for example, \cite[Theorem 3.1]{LO08} and \cite[Theorem 1]{O16}.

%%In this paper, we study the stability or the linearized stability of solutions $\vp_\pp(t,e)$ and $\phi_\pp(t,e)$.
Since the linearization equations of \x{xe} are Hill's equations with parameter $e$ \cite{MW66}, the linearized stability/instability of these periodic solutions $\vp_\pp(t,e)$ and $\phi_\pp(t,e)$ are related with the traces $\tau_\pp(e)$ of the corresponding Poincar\'e matrixes. For $e=0$, one has $\tau_\pp(0)=2$ because Eq. \x{x} is autonomous and $\vp_\pp(t)$ and $\phi_\pp(t)$ are parabolic. Hence the signs of $\tau'_\pp(0)=\f{\rd \tau_\pp(e)}{\rd e}|_{e=0}$, if they are nonzero, can yield the linearized stability or instability.
%%Consequently, these periodic solutions are linearized stable (respectively, unstable) for all $e>0$ small if $\tau'_\pp(0)<0$
%%(respectively, $\tau'_\pp(0)>0$).
As for even $(\pp)$-periodic solutions $\vp_\pp(t,e)$, a formula of $\tau'_\pp(0)$ has been obtained in \cite{ZCC18} and will be restated as \x{tau-e1} of this paper.

One of the main results of this paper is to derive the corresponding formula of $\tau'_\pp(0)$ for odd $(\pp)$-periodic solutions $\phi_\pp(t,e)$. See formula \x{dtau0} in \S \ref{criteria}. Note that formulas \x{dtau0} and \x{tau-e1} for $\tau'_\pp(0)$ are involved of nonconstant periodic solutions $\vp_\pp(t)$ and $\phi_\pp(t)$ of the autonomous equation \x{x}, which are not known explicitly.
%% and the perturbation term $F(x,t,e)$ of \x{xe}.

By applying these formulas to the elliptic Sitnikov problem $(S_e)$, we can obtain the following analytical results on the stability or instability for some families of symmetric periodic solutions.

    \bb{Theorem} \lb{M5-7}
For those frequencies $(\pp)=(2p,p)$ where $p\in \N$ is arbitrary, we have the following results.

{\rm (i)} For the odd $(2p,p)$-periodic solutions $\phi_{2p,p}(t,e)$, one has
    \(
    \tau'_{2p, p}(0)>0.
    \)
Consequently, for $e>0$ small, $\phi_{2p,p}(t,e)$ is hyperbolic and Lyapunov unstable.

{\rm (ii)} For the even $(2p,p)$-periodic solutions $\vp_{2p,p}(t,e)$, one has
    \(
    \tau'_{2p, p}(0)<0.
    \)
Consequently, for $e>0$ small, $\vp_{2p,p}(t,e)$ is elliptic and linearized stable.
    \end{Theorem}

It seems to us that these are the first analytical results on the stability or instability for the nonconstant symmetric periodic solutions of the elliptic Sitnikov problem $(S_e)$.

The organization of the paper is as follows. In \S \ref{Hill}, we will introduce some notions for Hill's equations. The linearization equations of autonomous equation \x{x} along symmetric periodic solutions will be discussed with the emphasis on the relation between the fundamental solutions of linearization equations and the solutions of Eq. \x{x} themselves. See Lemma \ref{psi12}. Moreover, a relation between the Poincar\'e matrixes and the period function of the periodic solutions of Eq. \x{x} will be found in Lemma \ref{hbn}. These results may be of independent interests. In \S \ref{criteria}, we will first give the bifurcation result on odd $(\pp)$-periodic solutions $\phi_\pp(t,e)$ of Eq. \x{xe}. See Theorem \ref{M1}. Then we will derive the formula of $\tau'_\pp(0)$ in Theorem \ref{M2}. Finally, in \S \ref{Sitnikov}, we will use the formulas of $\tau'_\pp(0)$ to analyze the elliptic Sitnikov problem $(S_e)$. The results of Theorem \ref{M5-7} will be proved in \S \ref{odd} and \S \ref{even}.

Note from Theorem \ref{M5-7} that we have only obtained analytical results for some families of symmetric periodic solutions with very specific frequencies $(\pp)=(2p,p)$, because we are dealing with nonconstant periodic solutions. In fact, it is found numerically and analytically in \cite{GNRR18, ZCC18} that the stability/instability depend on frequencies in a delicate way. As for the elliptic Sitnikov problem, we will prove in Theorems \ref{M4} and \ref{M6} that $\tau'_\pp(0)$ are always $0$ for both odd and even $(\pp)$-periodic solutions when frequencies $(\pp)$ satisfy $m/(2p)\not\in \N$.
%%For other frequencies, the stability/instabilty may depend on the second-order derivatives $\tau''_\pp(0)$ with respect to $e$. %%These delicate dependence phenomena of have been found As for the elliptic Sitnikov problem, besides the results in Theorem %%\ref{M5-7}, we have more. For example, for those  one has always
The remaining frequencies are $(\pp)=(2np, p)$, $n\ge 2$. For odd $(2np,p)$-periodic solutions, numerical simulation shows that $\tau'_{2np,p}(0)$ are always positive and $\phi_{2np,p}(t,e)$ will lead to instability. For even $(2np,p)$-periodic solutions, we will prove in Lemma \ref{rels} that the signs of $\tau'_{2np,p}(0)$ differ from that of the odd ones by a factor $(-1)^n$. Hence some even  solutions are linearized stable, while the others are unstable. These observations will be stated as a conjecture at the end of the paper.

\ifl
Let us mention some important ones, among a lot of existence results.

\buu In \cite[Theorem 3.1]{LO08}, Llibre and Ortega gave the following result. Let integers $p,\ m$ be in condition \x{mp}.
Then there exists $e_\pp\in(0,1]$ and a family of solutions of problem $(S_e)$, $\vp_\pp(t,e)$, $e\in [0,e_\pp)$ such that $\vp_\pp(t,e)$ is even, $\T$-periodic in $t$, $\vp_\pp(0,e)>0$, and has precisely $2p$ zeros in one period. Moreover,
a sharp estimate on the possible maximal eccentricity $e_\pp$ is also obtained there.

These are called the {\it even} periodic orbits of the $(\pp)$-type, because they have the shape like $\cos(p t/m)$. The original proof in \cite{LO08} is to use global continuation theory. However, if no estimate on $e_\pp$ is considered, these solutions can also obtained from the Implicit Function Theorem (IFT) with the starting periodic solutions $\vp_\pp(t) :=\vp_\pp(t,0)$ being nonconstant $\T$-periodic solutions of the circular Sitnikov problem $(S_0)$ of the $(\pp)$-type.

\buu With the same choice of $(\pp)$ as in \x{mp}, an appropriate translation of $\vp_\pp(t)$ in $t$ can lead to an odd, $\T$-periodic solution $\phi_\pp(t)$ of problem $(S_0)$ of the $(\pp)$-type. By the IFT again, one has then a smooth family of solutions of problem $(S_e)$, $\phi_\pp(t,e)$, $e\in [0,e_\pp^*)$ such that $\phi_\pp(t,e)$ is odd, $\T$-periodic in $t$, $\dot\phi_\pp(0,e)>0$, and has precisely $2p$ zeros in one period. For details, see Theorem \ref{M1}. These  are then called {\it odd} periodic solutions of the $(\pp)$-type.

\buu In \cite[Theorem 1]{O16}, Ortega proved a very interesting result, i.e. for any $m\in \N$, the above family $\phi_{1,m}(t,e)$ of odd periodic orbits of the $(1,m)$-type is {\it uniquely, globally} defined. That is, for any $e\in[0,1)$, such a $\T$-periodic solution $\phi_{1,m}(t,e)$ is existent and unique. Such a uniqueness is obtained from some property on solutions of the linearized equations satisfying the Dirichlet boundary conditions.

\fi

\section{Periodic Solutions and Linearization of Autonomous Equations} \setcounter{equation}{0} \lb{Hill}

\subsection{Periodic solutions of autonomous equations} We consider the autonomous equation \x{x} with the symmetries as before. By introducing
    \be \lb{Ex}
    E(x):=\int_0^x f(u) \du, \qq x\in \R,
    \ee
an even function such that $E(0)=0$ and $E(x)>0$ for $x\ne 0$, we know that solutions $x(t)$ of \x{x} satisfy
    \be \lb{ener}
    C_h: \qq \oh \dot x^2(t) +  E(x(t)) \equiv h,
    \ee
where $h\in [0, +\oo)$. For $h=0$, \x{ener} corresponds to the equilibrium $x(t)\equiv 0$. For
    \[
    0< h < E_{\max}:=\sup_{x\in \R} E(x),
    \]
$C_h$ consists of a nonconstant periodic orbit in the phase plane, whose minimal period is denoted by $T=T(h)>0$. We will not write down $T$ explicitly and refer to \cite{L91} for details.

Because of the symmetries of $f(x)$, we are interested in the following two classes of periodic solutions of Eq. \x{x}.

{\bf Odd periodic solutions:} For
    \be \lb{v}
    \eta\in \z(0,\eta_{\max}\y),\qq \eta_{\max}:=\sqrt{2 E_{\max}},
    \ee
let $x=\ss(t)=\ss(t,\eta)$ be the solution of \x{x} satisfying the initial value conditions
    \be \lb{ini}
    \z(x(0), \dot x(0)\y)=(0,\eta).
    \ee
Then $\ss(t)$ is a periodic solution of \x{x} of the minimal period
    \be \lb{Tv}
    T=T(h),\qq \mbox{where } h=\eta^2/2,
    \ee
with the following symmetries
    \be \lb{Os1}
    \ss(-t) \equiv -\ss(t) \andq \ss(t+T/2)\equiv - \ss(t).
    \ee
Moreover, $\ss(t)>0$ is strictly increasing on $(0,T/4)$.

{\bf Even periodic solutions:} For
    \[
    \xi\in \z(0,+\oo\y),
    \]
let $x=\cc(t)=\cc(t,\xi)$ be the solution of \x{x} satisfying the initial value conditions
    \be \lb{ini2}
    \z(x(0), \dot x(0)\y)=(\xi,0).
    \ee
Then $\cc(t)$ is a periodic solution of \x{x} of the minimal period
    \[
    T=T(h),\qq \mbox{where } h=E(\xi),
    \]
with the following symmetries
    \be \lb{Os2}
    \cc(-t)\equiv \cc(t) \andq \cc(t+T/2)\equiv - \cc(t).
    \ee
Moreover, $\cc(t)>0$ is strictly decreasing on $(0,T/4)$.

From \x{Os1} and \x{Os2}, one sees that
    \be \lb{Os12}
    \ss(T/2-t) \equiv \ss(t) \andq \cc(T/2-t) \equiv -\cc(t).
    \ee
The solutions $\ss(t)$ and $\cc(t)$ are also called $T/2$-anti-periodic. Like the sine and cosine, these solutions are related in the following way.

    \bb{Lemma} \lb{SC}
Suppose that $\eta$ and $\xi$ satisfy
    \be \lb{sc1}
     \eta^2/2=E(\xi)=:h.
    \ee
By setting $T=T(h)$, the odd and the even periodic solutions $\ss(t)=\ss(t,\eta)$ and $\cc(t)=\cc(t,\xi)$ are related via
    \be \lb{sc2}
    \ss(t+T/4) \equiv \cc(t) \andq \cc(t+T/4) \equiv -\ss(t).
    \ee
    \end{Lemma}

\subsection{Traces of Hill's equations} We need some general results for Hill's equations \cite{MW66}. Let $q: \R\to \R$ be a $T$-periodic locally Lebesgue integrable function and consider the Hill's equation
    \be \lb{he}
    \ddot y + q(t) y=0,\qq t\in \R.
    \ee
As usually, we use $y= \psi_i(t) =\psi_i(t,q)$, $i=1,2$ to denote the fundamental solutions of Eq. \x{he}, i.e. the solutions of \x{he} satisfying initial conditions $(\psi_1(0),\dot\psi_1(0)) = (1,0)$ and $(\psi_2(0),\dot\psi_2(0)) = (0,1)$ respectively.
The $T$-periodic Poincar\'e matrix of Eq. \x{he} is
    \[%be \lb{P-T}
    P=P_T=\matt{a}{b}{c}{d}:= \matt{\psi_1(T)}{\psi_2(T)}{\dot\psi_1(T)}{\dot\psi_2(T)}.
    \]%ee
The Liouville law for Eq. \x{he} asserts that
    \be \lb{Ll}
    \det P_T = a d - bc =+1.
    \ee
The trace of the $T$-Poincar\'e matrix $P_T$ is
    \[%be\lb{tau-T}
    \tau=\tau_T:= {\rm tr}(P_T)= a+d= \psi_1(T)+\dot\psi_2(T).
    \]%ee

Because of \x{Ll}, we know that (i) in case $|\tau|<2$, \x{he} is elliptic and is stable, (ii) in case $|\tau|>2$, \x{he} is hyperbolic and is unstable,  and (iii) the case $|\tau|=2$ corresponds to the parabolicity of Eq. \x{he} which can be either stable or unstable.

Being considered as functionals of potentials $q$, all of the above objects are Fr\'echet differentiable in $q\in L^1(\R/T\Z)$, the Lebesgue space endowed with the $L^1$ norm $\|\d\|_{L^1}$.

    \bb{Lemma} \lb{tau'0} {\rm (\cite[Lemma 2.2]{ZCC18})}
The Fr\'echet derivative of the trace $\tau: L^1(\R/T\Z)\to \R$ at $q$ is
    \be \lb{dtau}
    \f{\pa \tau}{\pa q}(h)=\int_0^T K(s) h(s) \ds\qqf h\in L^1(\R/T\Z).
    \ee
Here, by using the fundamental solutions $\psi_i(s)=\psi_i(s,q)$, %%the kernel is
    \be \lb{ker0}
    K(s):=-\psi_2(T) \psi^2_1(s) +\z(\psi_1(T)-\dot \psi_2(T)\y) \psi_1(s)\psi_2(s) +\dot \psi_1(T) \psi^2_2(s).
    \ee
    \end{Lemma}

\subsection{Linearization of autonomous equations} We consider a nonconstant $T$-periodic solution $x=\phi(t)$ of the autonomous equation \x{x}. Here $T$ is not necessarily the minimal period of $\phi(t)$. Then the linearization equation of \x{x} along the solution $\phi(t)$ is the Hill's equation \x{he}, where
    \be \lb{qt}
    q(t):= f'(\phi(t))
    \ee
is a $T$-periodic potential.

In the sequel, we consider
    \be \lb{pss}
    \phi(t):=\ss(t,\eta)\andq q(t):= f'(\ss(t,\eta)).
    \ee
Here $\ss(t,\eta)$ is an odd periodic solution of \x{x} of the minimal period $T$ as in \x{Tv}. Then one has the following important observations.

    \bb{Lemma} \lb{psi12}
Using the solutions $\ss(t,\eta)$ of initial value problems, the fundamental solutions $\psi_i(t)=\psi_i(t,q)$ of Eq. \x{he} are given by
    \bea\lb{psi10}
    \psi_1(t)\EQ %%\f{\dot \phi(t)}{\dot\phi(0)} =
    \f{1}{\eta}\pas{\f{\pa \ss}{\pa t}}{(t,\eta)} \andq
    \dot\psi_1(t) %%=\f{\ddot \phi(t)}{\dot\phi(0)}
    =-\f{f\z(\ss(t,\eta)\y)}{\eta},\\
    \lb{psi20}
    \psi_2(t)\EQ \pas{\f{\pa \ss}{\pa \eta}}{(t,\eta)}\andq \dot\psi_2(t)= \pas{\f{\pa^2 \ss}{\pa t\pa \eta}}{(t,\eta)}.
    \eea
    \end{Lemma}

\Proof Recall that $\ss(t,\eta)$ satisfies
    \bea \lb{phit}
    \EM\ddot \ss(t,\eta) +f\z(\ss(t,\eta)\y) = 0,\\
    \lb{ini-v}
    \EM \z(\ss(0,\eta),\dot \ss(0,\eta)\y)=(0,\eta).
    \eea
Differentiating \x{phit} with respect to $t$, we know that $y(t):=\pas{\f{\pa \ss}{\pa t}}{(t,\eta)}= \dot \ss(t,\eta)$ satisfies Eq. \x{he} and the initial values
    $$
    (y(0),\dot y(0))=\z(\dot \ss(0,\eta),\ddot \ss(0,\eta)\y)=\z(\dot \ss(0,\eta),-f\z(\ss(0,\eta)\y)\y)=\z(\eta,0\y) = \eta(1,0).
    $$
Hence we have
    \[
    \psi_1(t)\equiv{\dot \ss(t,\eta)}/\eta\andq \dot\psi_1(t)\equiv{\ddot \ss(t,\eta)}/\eta=-f\z(\ss(t,\eta)\y)/\eta,
    \]
the equalities in \x{psi10}.

On the other hand, by differentiating \x{phit} and \x{ini-v} with respect to $\eta$, we know that the variational equation for $y(t):= \pas{\f{\pa \ss}{\pa \eta}}{(t,\eta)}$ is just Eq. \x{he} and the initial values are $(y(0),\dot y(0))=(0,1)$. Thus $\psi_2(t) \equiv \pas{\f{\pa \ss}{\pa \eta}}{(t,\eta)}$. As a consequence,
    \[
    \dot\psi_2(t) \equiv \f{\pa}{\pa t}\z(\pas{\f{\pa \ss}{\pa \eta}}{(t,\eta)}\y)
    =\pas{\f{\pa^2 \ss}{\pa t\pa \eta}}{(t,\eta)}.
    \]
Thus we have the equalities in \x{psi20}.\qed

Since $f'(x)$ is even in $x$, it follows from \x{Os1} and \x{pss} that the minimal period of $q(t)$ is actually $T/2$. Because of this, we consider the Poincar\'e matrixes of Eq. \x{he} with different periods
    \[
    \hat P:=P_{T/2}\andq \hat{P}_n:=P_{n T/2},\q n\in \N.
    \]
Using the fundamental solutions $\psi_i(t)$, these are
    \[
    \hat P=\matt{\psi_1(T/2)}{\psi_2(T/2)}{\dot\psi_1(T/2)}{\dot\psi_2(T/2)}\andq \hat{P}_n=\matt{\psi_1(n T/2)}{\psi_2(n T/2)}{\dot\psi_1(n T/2)}{\dot\psi_2(n T/2)}.
    \]

    \bb{Lemma} \lb{PM}
By letting
    \be \lb{B1}
    \hat{b} :=\psi_2(T/2)\andq \hat b_n:=\psi_2(n T/2),
    \ee
one has
    \be \lb{Ptp}
    \hat P =\matt{-1}{\hat{b}}{0}{-1}\andq \hat P_n =\matt{(-1)^n}{\hat b_n}{0}{(-1)^n}, %% \andq P =\matt{1} {\psi_2(\T)} {0}{1}.
    \ee
and the constants $\hat{b}, \ \hat{b}_n$ are related via
    \be \lb{B2}
    \hat b_n=(-1)^{n+1}n\hat{b}.
    \ee
    \end{Lemma}

\Proof From \x{Os1} and their derivatives, one has
    \[
    \z(\ss(T/2),\dot \ss(T/2)\y)=\z(-\ss(0), -\dot \ss(0)\y)=\z(0, -\eta\y).
    \]
By \x{psi10}, we have
    $$
    \z(\psi_1(T/2),\dot\psi_1(T/2)\y)=\z(\dot \ss(T/2),-f(\ss(T/2))\y)/\eta=\z(-1,0\y),
    $$
i.e. the first column of $\hat P$ is $(-1,0)^\top$. Moreover, it follows from \x{Ll} that $\dot \psi_2(T/2)=-1$.  This gives the first result of \x{Ptp}.

For general $n\in \N$, one has then %%By letting the constants as in \x{B1}, one has from the first equality of \x{Ps}
    \[
     \hat{P}_n = \hat P^n= \matt{-1}{\hat{b}} {0} {-1}^n = \matt{(-1)^n}{{(-1)^{n+1}n \hat{b}}} {0} {(-1)^n}.
    \]
Hence  we have all equalities of the lemma.  \qed

Using the period function $T(h)$ of orbit $C_h$ of Eq. \x{x}, we have the following relation.

    \bb{Lemma} \lb{hbn}
Suppose that $T(h)$ is differentiable in $h$. Then
    \be \lb{B3}
    \hat{b}_n= (-1)^{n+1}n \f{\eta^2}{2} \z.\f{\rd T(h)}{\rd h}\y|_{\eta^2/2}= (-1)^{n+1}n h T'(h),
    \ee
where $h=\eta^2/2$ and $'=\f{\rd}{\rd h}$.
    \end{Lemma}

\Proof Since we are considering odd periodic solutions $\ss(t,\eta)$, we know from the second equality of \x{Os1} that
    \[
    \ss{(T(h)/2,\eta)} \equiv 0
    \]
for all $\eta$ as in \x{v}, where $h=\eta^2/2$ is as in \x{Tv}. Differentiating it with respect to $\eta$, we obtain
    \[
    \pas{\f{\pa \ss}{\pa t}}{(T(h)/2,\eta)} T'(h) \f{\eta}{2} + \pas{\f{\pa \ss}{\pa \eta}}{(T(h)/2,\eta)}=0.
    \]
By \x{psi10} and \x{psi20}, we have
    \[
    \pas{\f{\pa \ss}{\pa t}}{(T(h)/2,\eta)} =\eta \psi_1(T/2)=-\eta\andq \pas{\f{\pa \ss}{\pa \eta}}{(T(h)/2,\eta)}= \psi_2(T/2) = \hat b.
    \]
See the proof of Lemma \ref{PM}. Thus $\hat b=(\eta^2/2) T'(\eta^2/2)$. Combining with \x{B2}, we obtain result \x{B3} for general $n$.\qed

    \bb{Remark}\lb{nd20}
{\rm (i) From Lemmas \ref{PM} and \ref{hbn}, we have the following equivalence relations
    \be \lb{nd11}
\hat{b}\ne 0 \Longleftrightarrow \hat{b}_n\ne 0 \Longleftrightarrow  T'(h)\ne 0.
    \ee
One can notice that the former two conditions mean that $\phi(t)=S(t,\eta)$ is parabolic-unstable, while the last means that $\phi(t)$ is Lyapunov unstable because the periodic orbits inside a neighborhood of $C_h$ will have different periods.

(ii) For even periodic solutions $x=\cc(t)=\cc(t,\xi)$ of Eq. \x{x}, results analogous to those in Lemmas \ref{psi12}---\ref{hbn} have been deduced in \cite{ZCC18} in a similar way. }
    \end{Remark}

\section{A Stability Criterion for Odd Periodic Solutions} \setcounter{equation}{0} \lb{criteria}

\subsection{Bifurcations of odd periodic solutions} For $\eta>0$, we use $x=X(t,\eta,e)$ to denote the solution of problem \x{xe}-\x{ini}. In particular, when $e=0$, one has
    \be\lb{xs}
    X(t,\eta,0)\equiv \ss(t,\eta),
    \ee
the solution of problem \x{x}-\x{ini}.

Let $m\in \N$ and $p\in \N$. Suppose that there exists $h_\pp$ such that $C_{h_\pp}$ of \x{ener} is a periodic orbit of Eq. \x{x} of the minimal period $\T/p$, i.e.
    \be \lb{Tm}
    T(h_{\pp})= \T/p.
    \ee
Due to the autonomy and the symmetries of Eq. \x{x}, $C_{h_\pp}$ can be presented using either odd or even periodic solutions of Eq. \x{x}. In fact, by defining
    \be \lb{vpm}
    \phi_\pp(t):= \ss(t,\vpm), \qq \mbox{where }\vpm: = \sqrt{2 h_\pp},
    \ee
$\phi_\pp(t)$ is then an odd periodic solution of Eq. \x{x} of the minimal period $\T/p$. More symmetries of $\phi_\pp(t)$ can be found from \S 2.1. In particular, $\phi_\pp(t)$ is an $(\pp)$-periodic solution of \x{x} and satisfies
    \be \lb{ps0}
    \phi_\pp(t+m\pi/p)\equiv -\phi_\pp(t).
    \ee
This implies that $\phi_\pp(m \pi)=\ss(m\pi, \vpm)=0$, i.e. %%$\eta=\vpm$ is a solution of the following equation
    \be \lb{veq}
X(m\pi,\vpm,0) =0.
    \ee
See \x{xs}. As for the dependence of these solutions on $(\pp)$, one has $\phi_{mn,pn}(t) \equiv \phi_\pp(t)$ for any $n\in \N$.

A bifurcation result for odd $(\pp)$-periodic solutions of \x{xe} emanating from $\phi_\pp(t)$ is as follows.

    \bb{Theorem} \lb{M1}
Let $m, \ p$ and $h_\pp, \ \vpm$ be as above. Assume that
    \be \lb{nd90}
T'(h_\pp)\ne 0.
    \ee
Then there exist $e_\pp>0$ and a smooth function $E_\pp(e)$ of $e\in [0,e_\pp)$ such that
    \be \lb{Vs}
    E_\pp(0)=\vpm\andq X(m\pi,E_\pp(e),e)=0 \mbox{ for } e\in [0,e_\pp).
    \ee
Hence, for any $e\in [0,e_\pp)$,
    \be \lb{psie}
    \phi_\pp(t,e):= X(t,E_\pp(e),e)
    \ee
is an odd $(\pp)$-periodic solution of the non-autonomous equation \x{xe}, with the following symmetry
    \be \lb{p01}
    \phi_\pp(t+m\pi,e) \equiv -\phi_\pp(t,e).
    \ee
    \end{Theorem}

\Proof Let $\eta=\vpm$ be in Lemmas \ref{psi12}---\ref{hbn}. Then $T=\T/p$ and $m\pi= p \d T/2$. Thus
    \bea\lb{hbp}
    \pas{\f{\pa X}{\pa \eta}}{(m\pi,\vpm,0)}\EQ \pas{\f{\pa \ss}{\pa \eta}}{(m\pi,\vpm)}=\psi_2(m\pi)\qq \mbox{(by \x{psi20})}\nn\\
    \EQ \hat b_p \qq \mbox{(by \x{B1})} \nn\\
    \EQ (-1)^{p+1}p h_\pp T'(h_\pp)\qq \mbox{(by \x{B3})}\nn\\
    \NE 0\qq \mbox{(by \x{nd90})}.
    \eea
Combining with \x{veq}, the existence of the function $E_\pp(e)$ as in \x{Vs} follows immediately from the Implicit Function Theorem (IFT).

Since $F(x,t,e)$ is odd in $x$, the solution $\phi_\pp(t,e)$ of \x{psie} is obviously odd in $t$. Moreover, $\phi_\pp(t,e)$ satisfies \x{p01} and is $(\pp)$-periodic. \qed

    \bb{Remark} \lb{M11}
{\rm (i) As seen from \x{nd11} of Remark \ref{nd20}, the non-degeneracy condition \x{nd90} is equivalent to the instability of the $(\pp)$-periodic solution $\phi_\pp(t)$ of Eq. \x{x}, in the linearized sense and/or in the Lyapunov sense.

(ii) Note that $\phi_\pp(t,0) \equiv \phi_\pp(t)$ is $\T/p$-periodic. See \x{vpm}. Usually speaking, if $e>0$, the minimal period of $\phi_\pp(t,e)$ is $\T$, not $\T/p$.}
    \end{Remark}

\subsection{A stability criterion for odd periodic solutions} We consider the family $\phi_\pp(t,e)$ of odd $(\pp)$-periodic solutions of Eq. \x{xe} as in Theorem \ref{M1}.

For $e\in[0,e_\pp)$, the linearization equation of Eq. \x{xe} along $x=\phi_\pp(t,e)$ is the Hill's equation
     \be \lb{He}
    \ddot y + q(t,e)y=0, \qq q(t,e):=\pas{\f{\pa F}{\pa x}}{(\phi_\pp(t,e),t,e)}.
    \ee
Here the period is understood as $T=\T$. The corresponding trace is
    \be \lb{Tre}
\tau_\pp(e):={\psi_1(\T,e)}+{\dot\psi_2(\T,e)}.
    \ee
Here $\psi_i(t,e)$ are fundamental solutions of Eq. \x{He}. When $e=0$, we have
    \[%be \lb{pt}
    \phi_\pp(t,0)=\phi_\pp(t):= \ss(t,\vpm)\andq q(t,0)=q(t)= f'(\ss(t,\vpm)).
    \]%ee
See \x{pss}.

    \bb{Theorem}\lb{M2}
Let $\phi_\pp(t)$ be the odd $(\pp)$-periodic solution of Eq. \x{x} verifying condition \x{nd90}. Denote
    \be \lb{F23}
    F_{23}(t):=\pas{\f{\pa^2 F}{\pa t\pa e }}{(\phi_\pp(t),t,0)}.
    \ee
Then the derivative of the trace \x{Tre} at $e=0$ is
    \be \lb{dtau0}
    \tau'_\pp(0):= \pas{\f{\rd\tau_\pp(e)}{\rd e}}{e=0} =- p T'(h_\pp) \imt F_{23}(t)\dot\phi_\pp(t)\dt.
    \ee
Here $h_\pp=\vpm^2/2$ and $'=\f{\rd}{\rd h}$.
    \end{Theorem}

\Proof In order to apply Lemma \ref{tau'0}, we need to consider the $\T$-periodic Poincar\'e matrix $P$ of the linearization equation
    \[
    \ddot y+ q(t) y=0, \qq \mbox{where } q(t):=f'(\phi_\pp(t)).
    \]
Arguing as in the proof of \x{hbp}, by letting $T=\T/p$ in  Lemmas \ref{psi12}---\ref{hbn} and noticing that $\T=2p \d T/2$, we have
    \[
    \matt {\psi_1(2m\pi)} {\psi_2(2m\pi)} {\dot\psi_1(2m\pi)} {\dot\psi_2(2m\pi)} = \hat P_{2p}=\matt{1}{\hat b_{2p}}{0}{1},
    \]
where
    \be \lb{bmp}
    \hat b_{2p}= \psi_2(2m\pi)= -2 p h_\pp T'(h_\pp)=: b_\pp.
    \ee
See \x{B3} with $n=2p$. Thus the kernel of \x{ker0} is
    \[%be \lb{Kt}
    K(t) = -b_\pp \psi_1^2(t)= -\f{b_\pp}{\vpm^2} \dot\phi^2_\pp(t)\equiv p T'(h_\pp)\dot\phi^2_\pp(t).
    \]%ee

Denote
    \be \lb{phi}
    \Phi(t):=\pas{\f{\pa \phi_\pp(t,e)}{\pa e}}{(t,0)} \andq F_{13}(t)
    :=\pas{\f{\pa^2 F}{\pa e\pa x}}{(\phi_\pp(t),t,0)}.
    \ee
From \x{He}, we have
    \beaa
    h(t)\AND{:=}\pas{\f{\pa q}{\pa e}}{(t,0)} = \z.\f{\pa }{\pa e}\z(\pas{\f{\pa F}{\pa x}}{(\phi_\pp(t,e),t,e)}\y)\y|_{(t,0)} \\
    \EQ\pas{\f{\pa^2 F}{\pa x^2}}{(\phi_\pp(t),t,0)}\Phi(t) +F_{13}(t)\\
    \AND{=:} f''(\phi(t)) \Phi(t)+F_{13}(t).
    \eeaa
Here, for simplicity, $\phi(t):=\phi_\pp(t)$. From \x{dtau}, we obtain
    \bea \lb{tau0}
    \tau'_\pp(0)\EQ \imt K(t) h(t) dt
    %%=-b_\pp \inn{\psi_1^2}{f''(\phi) \Phi+F_{13}}
    = p T'(h_\pp)\imt\z(\Phi f''(\phi) + F_{13}\y)\dot\phi^2\dt.
    \eea

Since $\phi_\pp(t,e)$ is $\T$-periodic for any $e$, we know from the defining equality \x{phi} that $\Phi(t)$ is necessarily $\T$-periodic. Moreover, $\Phi(t)$ satisfies the variational equation
    \be \lb{phieq}
    \ddot{\Phi} + q(t)\Phi + F_{3}(t) =0,
    \ee
where
    \bea\lb{df3}
    F_{3}(t)\AND{:=}\pas{\f{\pa F}{\pa e}} {(\phi(t),t,0)},\nn\\
    \dot F_{3}(t)\EQ \f{\rd}{\dt}\z(\pas{\f{\pa F}{\pa e}} {(\phi(t),t,0)}\y)=
    \pas{\f{\pa^2 F}{\pa x \pa e}}{(\phi(t),t,0)}\dot \phi(t)
    + \pas{\f{\pa^2 F}{\pa t \pa e}}{(\phi(t),t,0)}\nn\\
    \EQ F_{13} (t) \dot \phi(t)+F_{23}(t).
    \eea

Recall that we have Eq. \x{phit} for $\phi(t)$ %%, Eq. \x{he} for $\dot\phi(t)\equiv \dot \phi(0) \psi_1(t)$
and Eq. \x{phieq} for $\Phi(t)$. From these we can obtain the following equality
    \be \lb{eqss}
    \f{\rd}{\dt}\z(\dot \Phi \ddot \phi- \ddot \Phi \dot \phi\y) =%%\Phi f''(\phi)\dot\phi^2 +\dot F_{3}\dot\phi\equiv
    \z(\Phi f''(\phi) + F_{13}\y)\dot\phi^2 +F_{23} \dot \phi.
    \ee
In fact, by using Eq. \x{phit} and Eq. \x{phieq}, one has
    \[
    \dot \Phi \ddot \phi- \ddot \Phi \dot \phi= -\dot \Phi f(\phi) + \Phi q \dot \phi + F_3 \dot \phi.
    \]
Thus the left-hand side of \x{eqss} is
    \beaa
    \EM -\f{\rd}{\dt}\z(\dot \Phi f(\phi)\y) + \f{\rd}{\dt}\z(\Phi q \dot \phi\y) + \f{\rd}{\dt}\z(F_3 \dot \phi\y)\\
    \EQ -\ddot \Phi f(\phi) - \dot \Phi f'(\phi)\dot\phi + \dot\Phi q \dot \phi+\Phi \dot q \dot \phi +\Phi q \ddot \phi+F_3 \ddot \phi +\dot  F_3 \dot \phi\\
    \EQ \z(\ddot \Phi+q \Phi +F_3\y) \ddot \phi +\z(-f'(\phi)+q \y) \dot \Phi \dot \phi +\Phi \dot q \dot \phi+\dot  F_3 \dot \phi \qq \mbox{(by \x{phit})}\\
    \EQ \Phi \dot q \dot \phi+\dot  F_3 \dot \phi\qq \mbox{(by \x{phieq} and \x{qt})}\\
    \EQ \Phi f''(\phi)\dot\phi^2 + F_{13}\dot\phi^2 +F_{23} \dot \phi\qq \mbox{(by \x{qt} and \x{df3})}.
    \eeaa

Finally, as $\Phi(t)$ and $\phi(t)$ are $\T$-periodic, by integrating \x{eqss} over $[0,\T]$, we obtain
    \[
    \imt \z(\Phi f''(\phi) + F_{13}\y)\dot\phi^2\dt +\imt F_{23} \dot \phi\dt =0.
    \]
Combining with \x{tau0}, we obtain the desired formula \x{dtau0}. \qed

Since $\tau_\pp(0)=2$, the role of formula \x{dtau0} is as follows.

    \bb{Corollary} \lb{M21}
{\rm (i)} If $\tau'_\pp(0)<0$, then $\phi_\pp(t,e)$ is elliptic and is linearized stable for $0<e\ll 1$.

{\rm (ii)} If $\tau'_\pp(0)>0$, then $\phi_\pp(t,e)$ is hyperbolic and is Lyapunov unstable for $0<e\ll 1$.
    \end{Corollary}

\subsection{A stability criterion for even periodic solutions, revisited} The bifurcations and linearized stability of even $(\pp)$-periodic solutions of Eq. \x{xe} have been done in \cite{ZCC18}. In the present notations, we restate the results in \cite{ZCC18} as follows. For $\xi>0$, we use $x=\X(t,\xi,e)$ to denote the solution of problem \x{xe}-\x{ini2}. Let $m\in \N$ and $p\in \N$ and the energy $h_\pp$ be as in \x{Tm}. By taking $\xi_\pp>0$ such that
    \[
    E(\xi_\pp)=h_\pp,
    \]
we know that
    \[%be \lb{vppm}
    \vp_\pp(t):= C(t,\xi_\pp)= \X(t,\xi_\pp,0)
    \]%ee
is an even $(\pp)$-periodic solution of \x{x} of the minimal period $T(h_\pp)=\T/p$. From Lemmas 2.5 and 2.6 of \cite{ZCC18}, under the same non-degeneracy condition \x{nd90}, i.e. $T'(h_\pp)\ne0$, one has from the IFT a smooth function $\Xi_\pp(e)$ of $e\in [0,\ul{e}_\pp)$ such that $\Xi_\pp(0)=\xi_\pp$ and
    \[
    \dot \X(m\pi,\Xi_\pp(e),e)\equiv 0.
    \]
Thus
    \[
    \vp_\pp(t,e) := \X(t,\Xi_\pp(e),e)
    \]
defines a family of even $(\pp)$-periodic solutions of Eq. \x{xe} which are emanated from $\vp_\pp(t)$. Moreover, $\vp_\pp(t,e)$ is also $m\pi$-anti-periodic as in \x{p01}.

Let $\ul{\tau}_\pp(e)$ be the trace of the $\T$-periodic Poincar\'e matrix of the linearization equation of \x{xe} along the solution $\vp_\pp(t,e)$. One has $\ul{\tau}_\pp(0)=2$ and the following formula.

    \bb{Theorem} \lb{M3} {\rm (\cite[Theorem 3.1]{ZCC18})}
With the notations above,
    \be \lb{tau-e1}
    \ul{\tau}'_\pp(0)= \z.\f{\rd\ul{\tau}_\pp(e)}{\rd e}\y|_{e=0} = - p T'(h_\pp) \imt\ul{F}_{23}(t) \dot\vp_\pp(t)\dt,
    \ee
where
    \be
    \lb{F23e}
    \ul{F}_{23}(t):=\pas{\f{\pa^2 F}{\pa t\pa e }}{(\vp_\pp(t),t,0)}.
    \ee
    \end{Theorem}

    \bb{Remark} \lb{M31}
{\rm For the case $m=1$, result \x{tau-e1} is proved in \cite{ZCC18}. See Formula (3.2) there. However, the coefficient there is expressed using $\dot {\ul{\psi}}_1(2\pi)$ and $f(\xi_{1,p})$, where $\ul{\psi}_1(t)$ is the first fundamental solution of the corresponding linearization equation. For general $m$, formula \x{tau-e1} can be deduced by a scaling of time. Moreover, arguing as in the deduction of \x{bmp}, the coefficient can be written in the present way. One can notice that the forms of formulas \x{dtau0} and \x{tau-e1} are the same.}
    \end{Remark}

\section{Stability Results for the Elliptic Sitnikov Problem} \setcounter{equation}{0} \lb{Sitnikov}

\subsection{Equations for the motions of the Sitnikov problems} After choosing the masses and the gravitational constant in an appropriate way, the governing equation for the motion of the infinitesimal mass in the elliptic Sitnikov problem $(S_e)$ is \cite{BLO94, LO08}
    \be \lb{se}
\ddot x +F(x,t,e) =0,\qq F(x,t,e):= \frac{x}{\z(x^2 + r^2(t,e) \y)^{3/2}}.
    \ee
Here $e\in [0, 1)$ is the eccentricity, and
    \be\lb{r}
    r(t,e) = r_0(1-e \cos u(t,e)), \qq r_0:=1/2,
    \ee
where, after some translation of time, $u=u(t,e)$ is the solution of the Kepler's equation
    \be \lb{Ke}
u - e\sin u =t.
    \ee
Note that the Kepler solution $u(t,e)$ is smooth in $(t,e)$ and satisfies
    \[%be \lb{ute}
    u(-t,e) \equiv - u(t,e)\andq u(t+2\pi,e) \equiv u(t,e) +2\pi.
    \]%ee
Consequently, $F(x,t,e)$ fulfills all requirements in \x{Sy1}. Moreover, when $e\in(0,1)$, the minimal period of $F(x,t,e)$ in $t$ is $2\pi$.

In particular, the circular Sitnikov problem $(S_0)$ is described by the autonomous equation
    \be \lb{s0}
    \ddot x + f(x)=0,\qq f(x):= \frac{x}{\z(x^2 + r_0^2\y)^{3/2}}.
    \ee
For Eq. \x{s0}, the energy $E(x)$ in \x{Ex} is
    \[
    E(x)=\int_0^x f(u) \du = 2- \f{1}{ \sqrt{x^2+r_0^2}}.
    \]
Solutions $x(t)$ of Eq. \x{s0} are on energy levels
    \be \lb{H}
    H(x,\dot x):=\f{1}{2} \dot x^2 - \f{1}{ \sqrt{x^2+r_0^2}}= h.
    \ee
Here the energy $h$ differs from that in \x{ener} by a constant $2$ and takes values from
     \(
     h\in[-2,+\oo).
     \)
For $h=-2$, \x{H} corresponds to the origin which is the equilibrium of \x{s0}. For $h\in(-2,0)$, \x{H} corresponds to periodic orbits of \x{s0} whose minimal period is denoted by $T(h)$. It is not difficult to verify that
    \[%be \lb{Ths}
    \lim_{h\to -2+} T(h) = 2\pi/\sqrt{8} \andq \lim_{h\to0-} T(h) =+\oo.
    \]%ee
Moreover, it is proved in \cite[Theorem C]{BLO94} that
    \be \lb{Th1}
    T'(h)=\f{\rd T(h)}{\rd h} >0 \qqf h\in (-2,0).
    \ee
Hence the origin is surrounded by a family of periodic orbits, whose minimal periods take values from $(2\pi/\sqrt{8},+\oo)$. For more facts on the dynamics of Eq. \x{s0}, see \cite{BLO94, LO08}.

To bifurcate the families $\phi_\pp(t,e)$ and $\vp_\pp(t,e)$ of $(\pp)$-periodic solutions of Eq. \x{se} which are respectively odd and even in $t$, the integers $m, \ p$ are required that $\T/p\in (2\pi/\sqrt{8},+\oo)$, i.e.
    \be \lb{mp1}
    1\le p \le \nu_m := [\sqrt{8} m],\qq m\in \N,
    \ee
because the non-degeneracy conditions \x{nd90} are ensured by \x{Th1}. Condition \x{mp1} is also used in \cite[\S 3]{LO08}. As before, we write $\phi_\pp(t,0)$ and $\vp_\pp(t,0)$ as $\phi_\pp(t)$ and $\vp_\pp(t)$ respectively. For these $(\pp)$-periodic solutions, it is convenient to call
    \[
    \vr:={p}/{m}
    \]
the rotation number. Condition \x{mp1} for $(\pp)$ is now equivalent to
    \be \lb{mp}
    \vr\in (0,\sqrt{8}) \cap \Q.
    \ee

\subsection{Analytical results for stability of odd periodic orbits}\lb{odd}
From the defining equalities \x{se}--\x{s0}, a direct computation can yield
    \be \lb{f23}
    \pas{\f{\pa^2 F}{\pa t\pa e}}{(x,t,0)} =\f{-3 x}{4\z(x^2+r^2_0\y)^{5/2}}\sin t.
    \ee
See also \cite[Formula (4.21)]{ZCC18}.

We first study the families $\phi_\pp(t,e)$ of odd $(\pp)$-periodic solutions of Eq. \x{se} for $m, \ p$ as in \x{mp1}. By \x{F23}, \x{dtau0} and \x{f23}, we have
    \[
    F_{23}(t)=\pas{\f{\pa^2 F}{\pa t\pa e}}{(\phi_\pp(t),t,0)} = \f{- 3 \phi_\pp(t)}{4\z(\phi_\pp^2(t)+r^2_0\y)^{5/2}}\sin t,
    \]
and
    \beaa
    \tau'_\pp(0)\EQ - p T'(h_\pp)\imt F_{23}(t)\dot \phi_\pp(t)\dt\\
    \EQ - \oq p T'(h_\pp)\imt \f{-3 \phi_\pp(t)\dot \phi_\pp(t)}{\z(\phi^2_\pp(t)+r^2_0\y)^{5/2}}\sin t\dt.
    \eeaa
Define
    \be \lb{G}
    G_\pp(t):={1}/{\z(\phi^2_\pp(t)+r^2_0\y)^{3/2}}.
    \ee
One has
    \[
    \dot G_\pp(t)= -{3 \phi_\pp(t)\dot \phi_\pp(t)}/{\z(\phi^2_\pp(t)+r^2_0\y)^{5/2}}.
    \]
Integrating by parts, we know that $\tau'_\pp(0)$ can be written as
    \be \lb{t11}
    \tau'_\pp(0) =\oq p T'(h_\pp)\imt G_\pp(t)\cos t\dt. %%\nn\\
    \ee
Such an observation was also used in \cite{ZCC18} for the study of even periodic solutions.

    \bb{Theorem} \lb{M4}
One has $\tau'_\pp(0)=0$ if $(\pp)$ satisfies \x{mp1} and
    \be \lb{pm0}
    \vr=\f{p}{m}\ne \f{1}{2}, \ \f{1}{4}, \ \f{1}{6},\ \dd
    \ee
In particular, $\tau'_\pp(0)=0$ if $m$ is odd and $1\le p\le \nu_m$, or $m$ is even and $m/2+1 \le p \le \nu_m$.
    \end{Theorem}

\Proof Let us notice from \x{ps0} and \x{G} that the minimal period of $G_\pp(t)$ is $m\pi/p$. Moreover, $G_\pp(t)$ is even in $t$. Hence one has the $m\pi/p$-periodic Fourier expansion
    \[
G_\pp(t) \equiv \sum_{n=0}^\oo a_n \cos\z(n \f{2p t}{m}\y)= \sum_{n=0}^\oo a_n \cos\z(2n p \f{t}{m}\y).
    \]
Let us write $\cos t$ as $\cos \z( m \f{t}{m}\y)$. By using the orthogonality of $\{\cos \z(n \f{t}{m}\y): n \in \Z^+\}$ in the space $L^2[0,\T]$, we know from \x{t11} that  $\tau'_\pp(0)=0$ if $(\pp)$ satisfies $m\ne 2 n p$ for all $n\in \N$, i.e. if $\vr$ satisfies \x{pm0}.
\qed

\ifl
As a corollary, one has the following results on the odd $\T$-periodic solutions.

    \bb{Corollary}\lb{M41}
One has $\tau'_\pp(0)=0$ if

\bu $m$ is odd and $1\le p\le \nu_m$, or

\bu $m$ is even and $m/2+1 \le p \le \nu_m$.
    \end{Corollary}
\fi

    \bb{Remark} \lb{M42}
{\rm From Theorem \ref{M4}, the signs of $\tau'_\pp(0)$ depend on the frequencies $(\pp)$ in a delicate way. For example, we have no information on the stability of odd $(\pp)$-periodic orbits $\phi_{m,p}(t,e)$ for any odd number $m$. This phenomenon was also observed for the families $\vp_\pp(t,e)$ of even periodic solutions of Eq. \x{xe} and Eq. \x{se}. See  \cite{ZCC18} and \cite{GNRR18}.}
    \end{Remark}

In contrast to case \x{pm0}, we have $m/(2p)=n \in \N$, i.e. $m=2p n$, or equivalently,
    \be \lb{pm09}
    \vr=\f{1}{2n},\qq n\in \N.
    \ee
In this case,
    \be \lb{phint}
    \phi_{2pn,p}(t)\equiv\phi_{2n,1}(t)=:\phi_n(t),
    \ee
which are the odd periodic solutions used by Ortega \cite{O16}. Note that $\phi_n(t)$ has the minimal period $T=\T/p=4n\pi$. More symmetries on $\phi_n(t)$ include
    \be \lb{sys}
    \z\{ \ba{l} \phi_n(-t) \equiv -\phi_n(t), \\
    \phi_n(t+2n\pi)\equiv -\phi_n(t),\\
    \phi_n(2n\pi-t) \equiv \phi_n(t), \\
    \phi_n(t)> 0 \q \mbox{for } t\in(0,2n\pi),\\
    \phi_n(t) \mbox{ is strictly increasing on $[0,n\pi]$.}
    \ea\y.
    \ee
Here the third equality of \x{sys} is deduced from \x{Os12}. Passing to the function
    \be \lb{Gnt}
    G_n(t):=1/\z(\phi_n^2(t)+r^2_0\y)^{3/2},
    \ee
one has
    \be \lb{sysn}
    \z\{ \ba{l} \mbox{$G_n(t)>0$ is even and  has the minimal period $2n\pi$,} \\
    G_n(2n\pi-t) \equiv G_n(t), \\
    \mbox{$G_n(t)$ is strictly decreasing on $[0,n\pi]$.}
    %%{ from $G_n(0)=8$ to $G_n(n\pi)= - h^3\in(0,8)$.}
    \ea\y.
    \ee
For the solution $\phi_n(t)$ as in \x{phint}, we can use the symmetries in \x{sysn} to obtain
    \beaa
    \int_0^{\T} G_n(t) \cos t\dt\EQ\int_0^{2p\d 2n\pi} G_n(t) \cos t\dt\\
    \EQ 2p \int_0^{2n\pi}G_n(t) \cos t\dt\nn\\
    \EQ 2 p \z(\int_0^{n\pi} G_n(t) \cos t\dt +\int_{n\pi}^{2n\pi} G_n(t) \cos t\dt\y)\\
    \EQ 4 p \int_0^{n\pi} G_n(t) \cos t\dt,
    \eeaa
because both $G_n(t)$ and $\cos t$ are symmetric with respect to $t=n\pi$. Combining with \x{Th1} and \x{t11}, we have the following results.

    \bb{Lemma}\lb{same}
For any $p, \ n\in \N$, we have
    \be \lb{t12}
    \tau'_{2pn, p}(0) = p^2 T'(h_{2n,1})A_n,
    \ee
where
    \be\lb{An}
    A_n:=\int_0^{n\pi} G_n(t) \cos t\dt= \oh \int_0^{2n\pi} G_n(t) \cos t\dt.
    \ee
In particular, $\tau'_{2pn, p}(0)$ and $A_n$ have the same sign for any $p\in \N$.
    \end{Lemma}

\ifl
    \bb{Theorem} \lb{M5}
For any $p\in \N$ and $m=2p$, i.e. $n=1$ in \x{pm09}, one has
    \[
    \tau'_{2p, p}(0)>0.
    \]
Consequently, for $e>0$ small, $\phi_{2p,p}(t,e)$ is hyperbolic and therefore is Lyapunov unstable.
    \end{Theorem}
\fi

Now we can complete the proof of Theorem \ref{M5-7} (i) for odd $(2p,p)$-periodic solutions $\phi_{2p,p}(t,e)$. The frequencies $(\pp)=(2p,p)$ correspond to the rotation number $\vr=\oh$. See \x{pm09}. Due to Lemma \ref{same}, we need only to prove that $A_1>0$. By \x{An}, one has $n=1$ and
    \bea \lb{A1}
    A_1\EQ \int_0^{\pi/2} G_1(t) \cos t\dt +\int_{\pi/2}^{\pi} G_1(t) \cos t\dt \nn\\
    \EQ \int_0^{\pi/2} G_1(t) \cos t\dt +\int_{\pi/2}^0 G_1(\pi-t) \cos (\pi-t)\rd(\pi-t)\nn\\
    \EQ \int_0^{\pi/2} \z(G_1(t)-G_1(\pi-t) \y) \cos t \dt. %\\
    \eea
From the last property of \x{sysn}, $G_1(t)$ is strictly decreasing on $[0,\pi]$. Hence \x{A1} implies that $A_1 >0$. \qed
%%%% Obtained on 2019.3.1.

\subsection{Analytical results for stability of even periodic orbits} \lb{even}

Let $m, \ p$ be as in \x{mp1}. We are now studying the family $\vp_\pp(t,e)$ of even $(\pp)$-periodic solutions of Eq. \x{se}. By \x{tau-e1}, \x{F23e} and \x{f23}, we have
    \bea \lb{t22}
    \ul{F}_{23}(t)\EQ \pas{\f{\pa^2 F}{\pa t\pa e}}{(\vp_\pp(t),t,0)} = \f{-3\vp_\pp(t)}{4\z(\vp_\pp^2(t)+r^2_0\y)^{5/2}}\sin t,\nn \\
    \ul{\tau}'_\pp(0)\EQ -p T'(h_\pp)\imt \ul{F}_{23}(t)\dot \vp_\pp(t)\dt\nn\\
    \EQ -\oq p T'(h_\pp)\imt \f{-3 \vp_\pp(t)\dot \vp_\pp(t)}{\z(\vp^2_\pp(t)+r^2_0\y)^{5/2}}\sin t\dt\nn\\
    \EQ \oq p T'(h_\pp)\imt \ul{G}_\pp(t)\cos t\dt,
    \eea
where
    \be \lb{Ge}
    \ul{G}_\pp(t):={1}/{\z(\vp^2_\pp(t)+r^2_0\y)^{3/2}}.
    \ee
Note that $\ul{G}_\pp(t)$ is even in $t$ and has the minimal period $m\pi/p$. The similar proof as in Theorem \ref{M4} can yield the following result.

    \bb{Theorem} \lb{M6}
One has $\ul{\tau}'_\pp(0)=0$ if $(\pp)$ satisfies \x{mp1} and \x{pm0}.
    \end{Theorem}

For the cases as in \x{pm09}, we have the following relation.

    \bb{Lemma} \lb{rels}
For any $p, \ n\in \N$, there holds
        \be \lb{rels1}
    \ul{\tau}'_{2pn,p}(0) = (-1)^n \tau'_{2pn,p}(0).
        \ee
    \end{Lemma}

\Proof We go back to formulas \x{t11} and \x{t22}, where $m=2pn$. Note that $\phi_\pp(t)$ and $\vp_\pp(t)$ have the same energy $h_{2pn,p}=h_{2n,1}$ and the same minimal period $T=\T/p=4n\pi$. Hence \x{sc1} is verified and the factors in \x{t11} and \x{t22} are the same. By \x{sc2}, one has
    \[
    \vp_\pp(t)\equiv \phi_\pp(t+n\pi).
    \]
By \x{Gnt} and \x{Ge}, we obtain the relation
    \[
    \ul{G}_\pp(t)\equiv G_\pp(t+n\pi)
    \]
Hence
    \beaa
    \imt \ul{G}_\pp(t) \cos t \dt \EQ \int_0^{4pn \pi}G_\pp(t+n\pi)\cos t\dt\\
    \EQ \int_{n\pi}^{n\pi+4pn \pi}G_\pp(t)\cos (t-n\pi)\dt\\
    \EQ (-1)^n \int_{n\pi}^{n\pi+4pn \pi}G_\pp(t)\cos t\dt\\
    \EQ (-1)^n \int_{0}^{4pn \pi}G_\pp(t)\cos t\dt,
    \eeaa
because $G_\pp(t)$ and $\cos t$ are $2n\pi$-periodic. Thus we have relation \x{rels1}. \qed

The stability result of Theorem \ref{M5-7} (ii) for even $(2p,p)$-periodic solutions $\vp_{2p,p}(t,e)$ follows immediately from Theorem \ref{M5-7} (i) and Lemma \ref{rels}. Hence the proof of Theorem \ref{M5-7} is complete.

\ifl
    \bb{Theorem} \lb{M7}
For any $p\in \N$ and $m=2p$, i.e. $n=1$ in \x{pm09}, one has
    \be \lb{ell2}
    \ul{\tau}'_{2p, p}(0)<0.
    \ee
Consequently, for $e>0$ small, $\vp_{2p,p}(t,e)$ is elliptic and therefore is linearized stable.
    \end{Theorem}

In fact, result \x{ell2} can be proved in a direct way. Arguing as in the deduction of \x{t12} and \x{A1}, we have
    \beaa
    \ul{\tau}'_{2p, p}(0)\EQ p^2 T'(h_{2,1})\int_0^{\pi} \ul{G}_1(t) \cos t\dt\\
    \EQ p^2 T'(h_{2,1})\int_0^{\pi/2} \z(\ul{G}_1(t)-\ul{G}_1(\pi-t) \y) \cos t \dt\\
    \LT 0,
    \eeaa
because, for the present case, the function
    \(
    \ul{G}_1(t) = {1}/{\z(\vp^2_{2,1}(t)+r^2_0\y)^{3/2}}
    \)
is strictly increasing on $[0,\pi]$.
\fi

\subsection{The numerical result and a conjecture} \lb{s44} For conservative systems like Hamiltonian systems, the stability of  periodic orbits is an important and a difficult problem \cite{SM71}. For the $N$-body problems and the related systems, one can refer to \cite{C10, C08, HLS14} for some different approaches to the stability of periodic orbits.

Going back to the Sitnikov problem, we know from Lemmas \ref{same} and \ref{rels} that, for any $n\ge 2$ and any $p\in \N$, the linearized stability/instability of $\phi_{2pn,p}(t)$ and $\vp_{2pn,p}(t)$ are determined by the sign of $A_n$. By \x{Gnt} and \x{An}, $A_n$ is only involved of the odd $(2n,1)$-periodic solution $\phi_n(t):=\phi_{2n,1}(t)$ of Eq. \x{s0}. It is easy to do the numerical simulation. With the choice of $1\le n \le 10$, we have the numerical results listed in Table \ref{tab1}.

\begin{table}
\centering
\caption{Numerical results for $\eta_n:=\eta_{2n,1}$, $h_n:=h_{2n,1}$ and $A_n$.} \lb{tab1}
\begin{tabular}{rrrrrrr}
$n$ & & $\eta_n$ & & $h_n$ & & $A_n$ \\
\midrule
1 & & $ 1.7192 $ & & $-0.5221$ & & $2.3179$  \\
2 & & $ 1.8319 $ & & $-0.3221$ & & $2.2194$  \\
3 & & $ 1.8735 $ & & $-0.2449$ & & $2.1843$  \\
4 & & $ 1.8965 $ & & $-0.2017$ & & $2.1615$  \\
5 & & $ 1.9112 $ & & $-0.1736$ & & $2.1479$  \\
6 & & $ 1.9216 $ & & $-0.1537$ & & $2.1380$  \\
7 & & $ 1.9294 $ & & $-0.1387$ & & $2.1293$  \\
8 & & $ 1.9355 $ & & $-0.1269$ & & $2.1227$  \\
9 & & $ 1.9404 $ & & $-0.1174$ & & $2.1174$  \\
10 & & $ 1.9445 $ & & $-0.1095$ & & $2.1131$  \\
\bottomrule[1pt]
\end{tabular}
\end{table}

\ifl

From \x{t11}, the signs of $\tau'_\pp(0)$ are only involved of the odd periodic solutions $\phi_\pp(t)$ of the circular Sitnikov problem \x{s0} we are considering. Thus it is easy to evaluate numerically. For $m=2,\ 4, \ 6, \ 8,\ 10, \ 12$, we list the numerical results in the table.

\bu the numerical result is consistent with the analytical results \x{pm0} and \x{A1}, and

\bu for all of those $(p,m)$ other than that in \x{pm0}, $\tau'_\pp(0)$ are always positive. It is then a very interesting question whether this is really  true.

\bu The problem is to prove the positiveness of
    \be\lb{An9}
    A_n=\int_0^{n\pi} \f{\cos t}{\z(\phi_n^2(t) +r_0^2\y)^{3/2}}\dt,\qq n=2,3,\dd
    \ee
where $\phi_n(t):=\phi_{2n,1}(t)$ are odd $(2n,1)$-periodic solutions considered in \cite{O16}, i.e. $\phi_n(t)$ is the unique odd $4n\pi$-periodic solution of the circular Sitnikov problem \x{s0} such that $\dot \phi_n(0)>0$ and $\phi_n(t)$ has the unique zero $t=2n\pi$ in the interval $(0,4n\pi)$.

From \cite{BLO94}, $\phi_n(t)$ can be expressed as elliptic functions of different kinds in an implicit way. Is this useful?
\fi

Note that the positiveness of $A_1$ in Table \ref{tab1} has already been proved in an analytical way. It is surprising that numerically, all of $A_n$, $n\ge 2$ are positive. Hence we have the following interesting problem.
\bigskip

\noindent
{\bf Conjecture} One has $A_n>0$ for all $n\ge 2$.
\bigskip

We end the paper with two remarks.

1. Once the conjecture is proved, we could conclude that (i) odd $(2np,p)$-periodic solutions $\phi_{2np,p}(t,e)$ are hyperbolic and Lyapunov unstable for $e>0$ small, (ii) even $(4np,p)$-periodic solutions $\vp_{4np,p}(t,e)$ are also  hyperbolic and Lyapunov unstable for $e>0$ small, and (iii) even $((4n-2)p,p)$-periodic solutions $\vp_{(4n-2)p,p}(t,e)$ are elliptic and linearized stable for $e>0$ small.

2. For the case $n=2$, arguing as in \x{A1}, we have from \x{An}
    \[
    A_2
    %%\EQ \int_0^{\pi/2} \z(\z(G_2(t)-G_2(\pi-t)\y) - \z(G_2(t+\pi)-G_2(2\pi-t) \y)\y) \cos t \dt\\
    =\int_0^{\pi/2} \z(G_2(t)-G_2(\pi-t)+G_2(2\pi-t) -G_2(\pi+t) \y) \cos t \dt.
    \]
The sign of $A_2$ is related with a certain kind of `convexity' of $G_2(t)$ on the interval $[0,2\pi]$. This is also true for general case $n\ge 3$.

\ifl
For general case $n\ge 2$, arguing as in \x{A1}, we have from \x{An}
    \bea\lb{An1}
    A_n\EQ \sum_{i=1}^n \int_{(i-1)\pi}^{i\pi} G_n(t) \cos t\dt\nn\\
    \EQ\sum_{i=1}^n \int_{0}^{\pi} G_n(t+(i-1)\pi) \cos (t+(i-1)\pi)\dt\nn\\
    \EQ \sum_{i=1}^n \int_{0}^{\pi} (-1)^{i-1} G_n(t+(i-1)\pi) \cos t\dt\nn\\
    \EQ \sum_{i=1}^n \int_{0}^{\pi/2} (-1)^{i-1} \z(G_n(t+(i-1)\pi) -G_n(i\pi-t)\y) \cos t\dt\nn\\
    \EQ \int_{0}^{\pi/2}\z( \sum_{i=1}^n (-1)^{i-1}\z(G_n(t+(i-1)\pi) -G_n(i\pi-t)\y)\y) \cos t\dt.
    \eea
It is crucial to prove that $A_2>0$. This is related with some `convexity' of $G_2(t)$. However, $G_2(t)$ cannot be convex in the whole interval $[0,2\pi]$.
\fi

%%%\section*{Acknowledgment}

\vfill \hfill \fbox{\small Ver. 1, 2019-04-26}

\end{document}